\documentclass{amsart}
\usepackage{amscd,amssymb}

\theoremstyle{plain}
\newtheorem{theorem}{Theorem}[section]
\newtheorem{proposition}[theorem]{Proposition}

\newtheorem{lemma}[theorem]{Lemma}

\begin{document}
\setcounter{section}{-1}

\title{Nice enumerations of $\omega$-categorical groups} 
\author{A.Ivanov and K.Majcher} 


\maketitle

\centerline{Institute of Mathematics, University of Wroclaw, Wroclaw, 
Poland 
\footnote{Correspondence: A.A.Ivanov, Institute of Mathematics, 
University of Wroclaw, pl.Grunwaldzki 2/4, 50-384, Wroclaw, Poland; 
E-mail: ivanov@math.uni.wroc.pl}
}
\centerline{Institute of Mathematics, University of Wroclaw, Wroclaw, 
Poland} 

\bigskip

\centerline{{\bf ABSTRACT}}

We give an example of an $\omega$-categorical group without 
AZ-enumerations. 
We build AZ-enumerations of some $\omega$-categorical 
central products of $\omega$ copies of a finite 2-step nilpotent group. 
\bigskip 



\section{Introduction} 
The following notion has been introduced by Alhbrandt and Ziegler in 
\cite{AZ} as a technical tool for quasifinite axiomatizability. 
An ordering $<$ of type $\omega$ of a countable structure $M$ is 
called a {\em nice enumeration} of $M$ if for any sequence 
$a_i \in M$, $i\in \omega$, there are $i,j\in \omega$ and some 
automorphism $\alpha \in Aut(M)$, such that 
$\alpha (a_i )=a_j$, and $\alpha(a)< a_j$ for all $a< a_i$.   
This notion has been applied in several places of model theory. 
The following question is central in the subject: 
\begin{quote} 
Is there an $\omega$-categorical structure without a nice enumeration ?
\end{quote} \parskip0pt 

In our paper we will study some version of nice enumerations 
which was introduced by Hrushovski in \cite{Hr}. 
An $\omega$-ordering $<$ of a countable structure $M$ is called 
an {\em AZ-enumeration} of $M$ if for each $n$ and any sequence 
$\bar{a}_i$, $i \in \omega$, of $n$-tuples of $M$ there are 
$i,j\in \omega$ and some order preserving elementary map 
$\alpha :M \rightarrow M$ such that 
$\alpha (\bar{a}_i )= \bar{a}_j$. 
In \cite{CHr} structures having AZ-enumerations are called 
{\em geometrically finite}. \parskip0pt 

It is easy to see that every AZ-enumeration is nice. 
Albert and Chowdhury have asked in \cite{AC} whether there is 
an $\omega$-categorical structure which is not geometrically finite.   
In particular they have asked if the random graph has an AZ-enumeration. 
In our paper we study AZ-enumerations in the case of 
$\omega$-categorical groups. 
We answer the question from \cite{AC} mentioned above by showing 
that the 2-step nilpotent group with quantifier elimination 
found in \cite{CSW} does not have AZ-enumerations. 
Using a similar idea we also prove that the random graph 
does not have AZ-enumerations. \parskip0pt 

In Section 2 we study the same questions for some 
$\omega$-categorical central products of $\omega$ copies of 
a finite 2-step nilpotent group (see \cite{A}). 
It is easy to see that these groups are reducts of smoothly 
approximable structures. 
By \cite{CHr} this implies that they have AZ-enumerations. 
We prove that the standard enumerations of these groups are 
already AZ-enumerations in some stronger sense.      
Although this theorem resembles some statements of Section 4.1 
from \cite{CHr}, our proof uses different ideas and moreover  
provides some additional information. \parskip0pt 

The research is supported by KBN grant 1 P03A 025 28. 

\section{Fraiss\'{e} limits without AZ-enumerations}

We start with the construction of a QE-group of nilpotency 
class 2 given in \cite{CSW}. 
Since the group is build as the Fraiss\'{e} limit of a class of 
finite groups, we give some standard preliminaries 
(see for example \cite{evans}). \parskip0pt 
   
Let $\mathcal{K}$ be a non-empty class of finite structures of 
some finite language $L$. 
We assume that $\mathcal{K}$ is closed under taking substructures 
(satisfies HP, the {\em hereditary property}), has the 
{\em joint embedding property} (JEP) and 
the {\em amalgamation property} (AP). 
The latter is defined as follows: for every pair of embeddings  
$e: A \rightarrow B$ and $f: A \rightarrow C$ with 
$A,B,C\in \mathcal{K}$ there are embeddings $g: B \rightarrow D$ 
and $h: C \rightarrow D$ with $D\in \mathcal{K}$ such that 
$g\cdot e =h\cdot f$. 
Fraiss\'{e} has proved that under these assumptions there is 
a countable locally finite $L$-structure $M$ 
(which is unique up to isomorphism) such that: 
\begin{quote}
(a) $\mathcal{K}$ is the {\em age} of $M$, i.e. the class of all 
finite substructures which can be embedded into $M$ and \parskip0pt 
 
(b) $M$ is {\em finitely homogeneous} (ultrahomogeneous), 
i.e. every isomorphism between finite substructures of $M$ 
extends to an automorphism of $M$.   
\end{quote}
The structure $M$ is called the {\em Fraiss\'{e} limit} of $\mathcal{K}$. 
\parskip0pt 

To define a $2$-step nilpotent, $\omega$-categorical group 
without AZ-enumerations we assume that $\mathcal{K}$ is 
the class of all finite groups of exponent four in which 
all involutions are central.  
By \cite{CSW} $\mathcal{K}$ satisfies the HP, the JEP and the AP. 
Let $\mathcal{G}$ be the Fraiss\'{e} limit of this class. 
Then $\mathcal{G}$ is nilpotent of class two. \parskip0pt 

We need the notions of free amalgamation and a-indecomposability 
in $\mathcal{K}$. 
Following \cite{CSW} we define them through the associated 
cathegory of {\em quadratic structures}. 
A quadratic structure is a structure $(U,V;Q)$ where $U$ and 
$V$ are vector spaces over the field $\mathbb{F}_2$ and 
$Q$ is a nondegenerate quadratic map from $U$ to $V$, i.e. 
$Q(x) \neq 0$ for all $x \neq 0$ and the function 
$\gamma(x,y)=Q(x)+Q(y)+Q(x+y)$ is an alternating 
bilinear map. 
By $\mathcal{Q}$ we denote the category of all 
quadratic structures with morphisms 
$(f,g):(U_1,V_1;Q_1) \rightarrow (U_2,V_2;Q_2)$
given by linear maps $f:U_1 \rightarrow U_2$, 
$g:V_1 \rightarrow V_2$ respecting the quadratic 
map: $g Q_1 =Q_2 f$. \parskip0pt 

For $G \in \mathcal{K}$ define $V(G):= \Omega (G)$, 
the subgroup of all involutions of $G$, and $U(G):= G/V(G)$. 
Let $Q_G: U(G)\rightarrow V(G)$ be the map induced by squaring in $G$. 
Then $QS(G)= (U(G),V(G);Q_G)$ is a quadratic structure 
and the associated map $\gamma (x,y)$ is the one induced by 
the commutation from $G/V(G) \times G/V(G) $ to $V(G)$. 
It is shown in Lemma 1 of \cite{CSW} that this gives 
a 1-1-correspondence between $\mathcal{K}$ and $\mathcal{Q}$ 
up to the equivalence of central extensions 
$1\rightarrow V(G)\rightarrow G\rightarrow U(G)\rightarrow 1$ 
with $G\in \mathcal{K}$. \parskip0pt    

We now consider the amalgamation process in $\mathcal{K}$. 
To any amalgamation diagram in $\mathcal{K}$, 
$G_0 \rightarrow G_1 ,G_2$ we associate the diagram 
$QS(G_0 )\rightarrow QS(G_1 ),QS(G_2 )$ of the corresponding 
quadratic structures and (straightforward) morphisms.   
Let $QS(G_i )=(U_i ,V_i ;Q_i )$, $i\le 2$. 
Let $U^*, V^*$ be the amalgamated direct sums 
$U_1 \bigoplus_{U_{0}} U_2$, $V_1 \bigoplus_{V_{0}} V_2$ 
in the category of vector spaces. 
We define the free amalgam of $QS(G_1 )$ and $QS(G_2 )$ 
as a quadratic structure $(U ,V ;Q)$ with $U=U^*$ and 
$V=V^* \bigoplus (U_1 /U_0 )\otimes (U_2 /U_0 )$ (see \cite{CSW}). 
The corresponding quadratic map $Q:U \rightarrow V$ 
is defined by first choosing splittings of $U_1$, $U_2$ 
as $U_0 \bigoplus U'_1$ and $U_0 \bigoplus U'_2$, 
respectively, identifying $U'_1$, $U'_2$ 
with $U_1 /U_0$, $U_2 /U_0$ and defining 
$$
Q(u_0 +u'_1 +u'_2 )= Q_0 (u_0) +Q_1 (u'_1) +Q_2 (u'_2) + 
\gamma_1 (u_0, u'_1 ) + \gamma_2 (u_0,u'_2 ) + (u'_1 \otimes u'_2 ).
$$
Note that $Q \upharpoonright U_i =Q_i $ and the corresponding 
$\gamma (u'_1, u'_2)$ is $u'_1 \otimes u'_2$. 
Since $u'_1 \otimes u'_2 =0$ only when one of the factors 
is zero, the nondegeneracy is immediate.
It is shown in \cite{CSW} that $(V,U;Q)$ is a pushout of 
the natural maps $QS(G_1 )$, $QS(G_2 )\rightarrow (V,U;Q)$ 
agreeing on $QS(G_0 )$. 
We call the quadratic structure $(V,U;Q)$ the {\em free amalgam} 
of $QS(G_1 )$, $QS(G_2 )$ over $QS(G_0 )$. 
Let $G$ be the group associated with $(V,U;Q)$ in $\mathcal{K}$. 
By Lemma 3 of \cite{CSW} there are embeddings 
$G_1, G_2 \rightarrow G$ with respect to which $G$ becomes 
an amalgam of $G_1$, $G_2$ over $G_0$ in $\mathcal{K}$. 
We call $G$ {\em the free amalgam} of $G_0 \rightarrow G_1 ,G_2$. 
\parskip0pt 

We call a group $H\in \mathcal{K}$ {\em a-indecomposable} 
if whenever $H$ embeds into the free amalgam of two structures 
over a third, the image of the embedding is contained in one 
of the two factors. 
It is proved in Section 3 of \cite{CSW} that there is 
a sequence of a-indecomposable groups 
$\{ G_d :d\in\omega\} \subseteq \mathcal{K}$ such that for 
any pair $d\not= d'$ the group $G_d$ is not embeddable into $G_{d'}$.  

\begin{theorem}
Let $\mathcal{G}$ be the Fraiss\'{e} limit of the class 
$\mathcal{K}$ of all finite 2-step nilpotent groups of 
exponent four such that all involutions are central. 
Then $\mathcal{G}$ does not have AZ-enumerations. 
\end{theorem}

{\em Proof.}
Let $\{ G_d :d\in\omega\}$ be an antichain of 
a-indecomposable groups as above. 
We define the rank $rk(G_d )$ as the minimal size of 
a generating set of $G_d$. 
Since $\mathcal{G}$ is $\omega$-categorical and all $G_d$ are 
embeddable into $\mathcal{G}$, we may assume that 
$2<rk(G_d )<rk(G_{d'})$ for every pair $0\le d<d'$. \parskip0pt 
 
Let $<$ be an ordering of $\mathcal{G}$. 
Suppose for a contradiction that $<$ defines 
an AZ-enumeration of $\mathcal{G}$. 
We will now define a sequence of triples $a_n <b_n <c_n$,  
$n\in \omega$, and a subsequence $G_{d_n}$, 
$n\in\omega\setminus \{ 0\}$, satisfying the following 
conditions. 
Let $a_n$ be the $<$-minimal element of $\mathcal{G}$ 
(thus $a_{0} =...=a_n =...$). 
For $n>0$ the elements $b_n$ and $c_n$ are chosen so that 
there is a subset $T_n$ consisting of some 
$x_1 ,x_2 ,...,x_{t_n}\le b_n$, 
such that the set $T_n \cup\{ c_n\}$ generates a subgroup 
isomorphic to $G_{d_n}$. 
We also demand that for each $i<n$ there is no subset 
$T \subseteq \{x: x \leq b_n \}$ such that 
$T\cup\{ c_n\}$ generates a subgroup isomorphic to $G_{d_i}$. 
\parskip0pt 

The triples $(a_n ,b_n ,c_n )$ are defined by induction. 
Let $a_0 =b_0 =c_0$. 
At step $n$ we take $b_n$ as the first element enumerated 
after $c_{n-1}$ such that the initial segment 
$\mathcal{G}_n =\{ x:x\leq b_n\}$ contains a set 
$T_n$ which together with some 
$c\in\mathcal{G}\setminus\langle\mathcal{G}_n \rangle$ 
generates a subgroup isomorphic to some 
$G_d \not\in\{ G_{d_1},...,G_{d_{n-1}}\}$. 
Let $d_n$ be the minimal number $d$ with this condition.  
To define $c_n$ consider a group $U_n$ which 
is isomorphic to the free amalgam of $G_{d_n}$ and 
$\langle\mathcal{G}_n\rangle$ over $\langle T_n\rangle$ by 
an isomorphism fixing $\mathcal{G}_n$ pointwise. 
Since $\mathcal{G}$ is the Fraiss\'{e} limit of $\mathcal{K}$ 
we see that $U_n$ can be chosen as $\langle \mathcal{G}_n ,c \rangle$ 
for an appropriate $c\in \mathcal{G}$. 
Let $c_n$ be the element of $\mathcal{G}$ with the minimal number 
with respect the condition that  
$\langle \mathcal{G}_n , c_n \rangle$ is isomorphic with 
$\langle \mathcal{G}_n ,c \rangle$ over $\mathcal{G}_n$ 
under an isomorphism taking $c_n$ to $c$. \parskip0pt 

{\em Claim.} 
There are no $i<n$ and a subset $T\subseteq \mathcal{G}_n$ 
such that $T\cup\{ c_n\}$ generates a subgroup isomorphic to $G_{d_i}$. 

Suppose that such $T$ exists. 
This defines a copy of $G_{d_i}$ in the free amalgam of $G_{d_n}$ 
and $\langle\mathcal{G}_n \rangle$ over $\langle T_n \rangle$. 
By a-indecomposability either 
$\langle T\cup \{ c_n\} \rangle \subseteq G_{d_n}$ or 
$\langle T\cup \{ c_n\} \rangle \subseteq \langle \mathcal{G}_n \rangle$.
The first case is impossible because there is no embedding 
of $G_{d_i}$ into $G_{d_n}$. 
The second condition contradicts the assumption that 
$c_n \not\in \langle \mathcal{G}_n\rangle$.
\parskip0pt 

To finish the proof of the theorem assume that   
$\rho :M \rightarrow M$ is an order preserving elementary map 
taking $(a_i,b_i,c_i)$ to $(a_j,b_j,c_j)$ for some $0<i<j$. 
Since $\langle T_i \cup\{ c_i \}\rangle$ is isomorphic to 
$G_{d_i}$, there is a subset $T\subseteq \mathcal{G}_j$ such 
that $\langle T\cup\{ c_j\}\rangle$ is isomorphic to $G_{d_i}$ 
(for example let $T=\rho (T_i)$).  
This contradicts the definition of triples $(a_n ,b_n ,c_n)$, 
$n\in\omega$. 
${\square}$ 

\bigskip 

We finish this section by a similar argument applied 
to graphs. 
Although it does not concern $\omega$-categorical groups, 
we have decided to include it into the paper.   
Besides the fact that this argument is very similar, 
it answers a question from \cite{AC}, 
which was somehow distinguished in that paper.    

Let $\mathcal{K}_0$ be the class of all finite graphs. 
The Fraiss\'{e} limit $(\Gamma ,R)$ of $\mathcal{K}_0$ 
is called the {\em random graph.}
 
\begin{proposition}
The random graph does not have AZ-enumerations. 
\end{proposition}

{\em Proof.}
Assume for a contradiction that there is an ordering $<$ 
which defines an AZ-numeration of the random graph $(\Gamma ,R)$. 
We define an infinite sequence of triples $a_n <b_n <c_n$, 
$n\in \omega\setminus \{ 0\}$, satisfying the following conditions. 
All $a_n$ always denote the (same) $<$-minimal element of $\Gamma$. 
For $n>3$ the elements $b_n$ and $c_n$ are chosen so 
that there are $x_1, x_2,..., x_{n}\le b_n$ such that 
$x_1,x_2,...,x_{n}$ form an $n$-cycle, i.e.  
$(x_i ,x_j )\in R\leftrightarrow (|j-i|=1)\vee (\{ i,j\}=\{ 1,n\} )$, 
and $(x_i ,c_n )\in R$ for all $i\le n$. 
On the other hand we demand that for each $i$ with $2<i<n$, 
any $y_1, y_2 ,...,y_{i}\le b_n$ does not form an $i$-cycle 
$R$-connected with $c_n$ as above.  \parskip0pt 

The triples $(a_n ,b_n ,c_n )$, $n\in\omega$, 
can be defined by induction. 
At step $n$ we take $b_n$ as the first element enumerated after 
$c_{n-1}$ such that there are $x_1 ,x_2 ,...,x_{n}\le b_n$ forming 
an $n$-cycle: $(x_i , x_j)\in R$ if only if $|i-j|=1$ 
or $\{ i,j\} =\{ 1,n\}$. 
To find $c_n$ consider a finite graph $G$ consisting 
of vertices $x_1 ,x_2 ,...,x_{n}$ and some $c$ with 
$(x_i ,c)\in R$ for $i\le n$. 
Let $B_n$ be the free amalgam of $\Gamma_n =\{ x:x\le b_n\}$ 
and $G$ over $\{ x_1 ,...,x_n\}$. 
Thus $B_n = \Gamma_n \cup\{ c\}$ where $c$ is not adjacent 
to any element of $\Gamma_n \setminus \{ x_1 ,...,x_{n}\}$. 
Since $\Gamma$ is homogeneous, the element $c$ can be found in $\Gamma$. 
Let $c_n$ be such an element $c\in\Gamma$ with the minimal 
number with respect to $<$. 
It is clear that for any $i<n$ there is no subset $T\subseteq \Gamma_n$ 
which forms an $i$-cycle $R$-connected with $c_n$. \parskip0pt 

Let $\rho : M \rightarrow M$ be an order preserving elementary map 
taking $(a_i,b_i,c_i)$ to $(a_j,b_j,c_j)$ for some $i<j$ with $3<i$. 
Since there is an $i$-cycle $x_1, x_2,..., x_{i}\le b_i$ such that 
$(x_{l},c_i)\in R$ for all $l\le i$, there is an $i$-cycle 
$T\subseteq \Gamma_j$ 
such that $T\cup \{ c_j\}$ is isomorphic to the structure 
defined on $\{ x_1 ,...,x_i ,c_i\}$. 
This contradicts the definition of triples $(a_n ,b_n ,c_n)$, 
$n\in\omega$. ${\square}$

\section{Some nilpotent groups with AZ-enumerations}
  
The following construction has been studied by Apps in \cite{A}. 
Suppose that $G$ is a group, $K$ is a subgroup of $Z(G)$, 
and $A$ is some indexing set of cardinality $\lambda$. 
We define $G(A;K)$, the (central) product of $\lambda$ 
copies of $G$ amalgamated over $K$, as follows.  
We denote by $G^A$ and $K^A$ the direct product of $\lambda$ 
copies of $G$ and $K$ respectively, indexed by $A$. 
Let $(K^A)^0 =\{ \gamma\in K^A :\Pi_{i\in A}\pi_{i}(\gamma ) =1\}$, 
where $\pi_{i} :K^A \rightarrow K$ is the projection map 
corresponding to $i\in A$, and let $G(A;K):= G^A/(K^A)^0$. 
We write $G(n;K)$ and $G(\omega;K)$ for $G(A;K)$ 
when $|A|=n$ and $\omega$ respectively. 
Note that if $i\in A$, then $G_{i}$, the $i$-th component 
of the direct product $G^A$, embeds into $G(A;K)$ over 
the identification map $K\rightarrow K^A /(K^A)^0$. 
The following theorem has been proved in \cite{A} (Theorem A). 
\begin{quote}
Let $G$ be finite, class 2 nilpotent group, and 
$K$ be a subgroup of $G$ such that $G' \leq K \leq Z(G)$. 
Then $G(\omega;K)$ is $\omega$-categorical.   
\end{quote}

We will improve this theorem by the statement that 
$G(\omega ;K)$ is a reduct of a smoothly approximable structure. 
This implies by \cite{CHr} that $G(\omega ;K)$ has 
an AZ-enumeration. \\ 
{\em Notation.}
Let $G$ and $K$ be as in the theorem. 
Denote $\Gamma =G(\omega;K)$, $\Gamma_n =G(n;K)$. 
Let $ \pi :G^{\omega} \rightarrow \Gamma$ be the quotient map. 
For $i\in \omega$, let $G_i$ be the $\pi$-image in $\Gamma$ 
of the $i$-th component of $G^{\omega}$ 
(which is also denoted by $G_i$). 
We have $G_i \cong G$ for each $i\in\omega$, and 
$\langle G_i: 0\le i<n\rangle$ is naturally isomorphic to $\Gamma_n$ 
(by the fact that $(K^{\omega})^{0} \cap G^{n} =(K^{n})^{0}$). 
We therefore view $\Gamma_n$ as a subgroup of $\Gamma$. 
\parskip0pt 

We now introduce a countable subgroup $\Omega <Aut(\Gamma )$ 
generated by some family of automorphisms studied in \cite{A}.
Let $\sigma$ be a finitary permutation of $\omega$. 
Then $\sigma$ induces an automorphism $\hat{\sigma}$ of 
$G^{\omega}$ given by 
$\hat{\sigma} (g_0,g_1,...) =(g_{\sigma (0)} ,g_{\sigma (1)},...)$. 
It is easy to see that $\hat{\sigma} (K^{\omega})^{0} =(K^{\omega})^{0}$. 
Thus $\hat{\sigma}$ can be considered as an automorphism of $G(\omega; K)$ 
such that $\hat{\sigma}(G_n)=G_{\sigma (n)}$ for each $n\in \omega$.
\parskip0pt 

Another kind of our automorphisms is defined as follows.  
Let $M= m+2$, where $m$ is the exponent of $G$. 
Let $\alpha_i :G\rightarrow G^M$ be given by 
$\alpha_i (g)=(g,...,|_i ,...,g)$, the $M$-tuple 
whose $i$-th entry is 1, and whose other entries are $g$. 
Define $\alpha: G^M \rightarrow G^M$ by 
$\alpha (g_0 ,...,g_{M-1} )=$ 
$\alpha_1 (g_0 )\cdot\cdot\cdot\alpha_M (g_{M-1} )$. 
Let $\beta =\pi \alpha: G^M \rightarrow \Gamma_M$. 
The following lemma has been proved in \cite{A} (Lemma 2.1). 

\begin{lemma} \label{lad} 
The map $\beta: G^M \rightarrow \Gamma_M$ is a homomorphism, 
and it induces an endomorphism $\beta^*$ of $\Gamma_M$ 
(i.e. $(K^{M})^{0} <ker \beta$). 
Moreover, $\beta^*$ is a self-inverse automorphism of 
$\Gamma_M$ which fixes every element of $K$
\end{lemma} 

It is worth noting that $\beta^* (g,1,...,1) =(1,g,...,g)$. 
We can consider $\beta^*$ as an automorphism of $\Gamma$ by 
defining its action trivially for entries with indexes 
greater than $M-1$. 
To see this it suffices to note that by Lemma \ref{lad} 
the kernel of the map $G^{\omega}\rightarrow \Gamma$ 
corresponding to this extension is contained in $(K^{\omega})^{0}$. 

\begin{lemma} \label{2type}
Let $i_0 ,j_0 \in \omega$ and $I$ be a finite subset of $\omega$ 
such that $\{ i_0 ,j_0\}\cap I=\emptyset$ and the exponent of $G$ 
divides $|I|$.
Then there is an automorphism $\alpha_{I,i_0 ,j_0}\in Aut(\Gamma )$ 
such that for any $(g_0 ,g_1 ,...,g_i ,...)$ from $G^{\omega}$ 
with $g_i =1$ for $i\in I\cup \{ j_0\}$, and $g_{i_0} \neq 1$,   
the automorphism $\alpha_{I,i_0 ,j_0}$ sends 
$(g_i: i< \omega)$ to $(g'_i: i< \omega)$, where $g'_i=g_i$ 
for $i \not\in I$ and $g'_i =g_{i_0}$ otherwise.
\end{lemma}

{\em Proof.}
The automorphism $\alpha_{I,i_0 ,j_0}$ can be chosen as 
a composition of automorphisms of the form $\hat{\sigma}$ 
for $\sigma\in S_{fin}(\omega )$ and automorphisms 
as in Lemma \ref{lad}. 
$\square$

\bigskip

Let $\Omega$ be the subgroup of $Aut(\Gamma )$ generated by all   
automorphisms as in Lemma \ref{2type} and all automorphisms 
of the form $\hat{\sigma}$ for $\sigma\in S_{fin}(\omega )$. 
\parskip0pt 

We start our study of $\Gamma$ with the observation that $\Gamma$ 
is a reduct of a smoothly approximable structure. 
We remind the reader that a structure $M$ is {\em smoothly approximable} 
if it is $\omega$-categorical and every finite subset of $M$ is 
contained in a finite substructure $N$ such that all $0$-definable 
relations on $M$ induce $0$-definable relations on $N$ and 
any two enumerations $\bar{a}$ and $\bar{b}$ of $N$ have 
the same type in $N$ if and only if they have the same type in $M$.   

\begin{proposition}
Let $G$ be a finite nilpotent group of class 2, and $K$ be 
a subgroup of $G$ such that $G' \leq K \leq Z(G)$. 
Then the constant expansion of $G(\omega ;K)$ by all elements 
of $K$ is smoothly approximable.
\end{proposition}

{\em Proof.} 
Consider all subgroups $\Gamma_n <\Gamma =G(\omega ;K)$, 
$n\in \omega$, realizable on the corresponding indexes $0,...,n-1$.  
We claim that $(\Gamma ,a)_{a\in K}$ is approximated 
by all $(\Gamma_n ,a)_{a\in K}$, $n\in\omega$. 
To see this it suffices to notice that every automorphism 
of $\Gamma_n$ fixing $K$ pointwise extends to an automorphism 
of $\Gamma$. 
Since $\Gamma$ is the central product of 
$G(\{ i\in\omega :i<n\};K)$ and $G(\{ i\in\omega :i>n-1\};K)$ 
amalgamated over $K$, we can extend an automorphism 
$\phi\in Aut(\Gamma_n /K)$ to $\Gamma$ 
trivially on $G(\{ i\in\omega :i>n-1\};K)$. 
$\square$

\bigskip 

We now build an explicit AZ-enumeration of $G(\omega ;K)$. 
In fact this is the corresponding version of the 
{\em standard ordering} of a basic linear geometry 
defined in Section 4.1 of \cite{CHr}. 
In our context this construction provides an AZ-enumeration 
with some additional properties. 
To formulate them consider a subgroup $\mathcal{H}$ of the group 
of all automorphisms of a structure $M$. 
The closure of $\mathcal{H}$ in the space $M^M$ of all functions 
$M\rightarrow M$ consists of some embeddings of $M$ into $M$. 
Since every element of $\mathcal{H}$ is an elementary map, 
these embeddings are elementary too. 
We say that an ordering $<$ of the structure $M$ is 
an {\em AZ-enumeration with respect to} $\mathcal{H}$       
if for any $n$ and any sequence $\bar{a}_i$, $i\in \omega$, 
of $n$-tuples from $M$ there are $i\not= j$ and some 
order preserving map from the closure of $\mathcal{H}$ 
which maps $\bar{a}_i$ to $\bar{a}_j$. \parskip0pt  

Let us define an AZ-enumeration of $G(\omega ;K)$. 
First we enumerate the group $G^{\omega}$.
Fix an ordering of $G$: $g_0$,...,$g_{m-1}$, where $g_0 =1$. 
Then we order $G^{\omega}$ by the reverse lexicographic 
ordering: $(a_1,a_2,...)< (b_1,b_2,...)$ if there is $j\in\omega$ 
such that $a_j<b_j$ and $a_i=b_i$ for all $i>j$. \parskip0pt 

We now construct an enumeration $\{ v_i: i< \omega\}$ 
of the group $\Gamma =G^{\omega}/(K^{\omega})^0$ by induction.
Suppose that $v_0,v_1,...v_{n-1}$ are already defined. 
Then let $v_n$ be the $(K^{\omega})^0$-coset 
having a representant which is minimal (with respect 
to the ordering above) in $G^{\omega}$ among 
sequences not representing $v_0 ,...,v_{n-1}$.  

\begin{theorem} \label{main2}
The ordering of the group $G(\omega ;K)$ defined as 
above is an AZ-enumeration with respect to $\Omega$. 
\end{theorem} 

Although this theorem corresponds to Lemma 4.1.6 from 
\cite{CHr}, our proof is based on some different tricks.   
When we construct a required order-preserving elementary 
map we explicitly define an approximating sequence 
from $\Omega$ guaranteeing that the map belongs to 
the closure of $\Omega$. 
It is possible that some special analysis of definable 
subsets of $G(\omega ;K)$ can be applied in this space 
instead of approximating sequences. 
However we think that our approach is more direct 
and elegant.\parskip0pt      

We start with some preliminaries. 
The following lemma belongs to G.Higman 
(see Section 4.1 of \cite{CHr}). 

\begin{lemma} \label{H}
Let $\Sigma$ be a finite set. 
Define a partial ordering on the set $\Sigma^{*}$ of 
$\Sigma$-words by: $w_1 \leq w_2$ if $w_1$ is a subword of $w_2$, 
i.e. after deleting some members of $w_2$ we are left with $w_1$. 
Then $( \Sigma^{*}, \leq )$ is a partial well ordering: for 
every sequence $\{ w_i :i<\omega\}$ from $\Sigma^{*}$, 
there are $i<j< \omega$ such that $w_i \leq w_j$. 
\end{lemma}

We now improve this lemma as follows. 
Consider again the set $\Sigma^{*}$ of all finite words over $\Sigma$. 
We say that a word $(a_i: i \leq n)$ is $*$-embedded into a word  
$(b_i: i \leq m)$ if there is an order preserving injection 
$f:\{ 1,...,n\} \rightarrow \{ 1,...,m\}$ such that $b_{f(i)} =a_i$ and 
$$
(\forall i\le m)(\exists j\le n)((i\le f(j))\wedge (b_i =b_{f(j)})).
$$  
It is easy to see that the following relation is 
a partial ordering on $\Sigma^*$: $w_1 \leq^{*} w_2$ if $w_1$ 
is $*$-embedded into $w_2$. 

\begin{lemma} \label{MH}
Let $\Sigma$ be a finite set. 
Then $( \Sigma^{*}, \leq^{*})$ is a partial well ordering.
\end{lemma}

{\em Proof.}
Let $A=\{ w_i: i\in\omega \} \subseteq \Sigma^*$.
We assume that all words in $A$ are composed from 
the same letters: $\sigma_1 ,...,\sigma_k$. 
Moreover we may also assume that for all $i\in \omega$ 
and $s<t \leq k$, the last appearance of $\sigma_s$ in $w_i$ 
is before the last appearance of $\sigma_t$ in $w_i$.
We thus view each word $w_i$ as 
$w_{i1}w_{i2}...w_{ik}\sigma_k$, where for $l>1$ the subword $w_{il}$ 
is of the form $\sigma_{l-1} \sigma_{r_1}\sigma_{r_2}...\sigma_{r_s}$ 
with $\sigma_{r_t} \in\{\sigma_l ,\sigma_{l+1},..., \sigma_k \}$, 
$t\le s$. 
It is enough to prove, that there are $i<j$, 
and an order preserving embedding 
$f:\{ 1,...,|w_i |\}\rightarrow \{ 1,...,|w_j |\}$, 
which sends $w_{il}$ to $w_{jl}$ for all $l<k$. 
For each $i \in\omega$ let $o_i$ be $max(|w_{i1}|,...,|w_{ik}|)$.   
In order to apply Lemma \ref{H} we will code up $w_i$ in some new alphabet. 
Let $\hat{\Sigma}$ be the alphabet of all $k$-tuples fom $G \cup \{x \}$. 
We associate to $w_i$ the word 
$\bar{\tau}^{(i)}=\tau_1 \tau_2 ...\tau_{o_i}$, 
where every $\tau_t$ is the sequence of $t$-th letters appearing  
in the corresponding $w_{il}$ (when $|w_{il}|<t$ the corresponding 
place in $\tau_t$ is signed by $x$). 
By Higman's lemma, there are $i<j<\omega$ and an embedding 
$f':\bar{\tau}^{(i)} \rightarrow \bar{\tau}^{(j)}$.   
Then $f'$ induces an $*$-embedding 
$f:\{ 1,...,|w_i |\}\rightarrow \{ 1,...,|w_j |\}$ of 
$w_i$ to $w_j$.  
To see this put $f(x)=y$ if there are $q \leq k$ and 
$r \le o_i$ such that $x= |w_{i1}|+|w_{i2}|+...+|w_{iq}|+r$,
$y= |w_{j1}|+|w_{j2}|+...+|w_{jq}| +f'(r)$. 
The rest is obvious. 
${\square}$

\bigskip 

{\em Proof of Theorem \ref{main2}.}
Let $\bar{a}_k \in (\Gamma_{\omega})^n$, $k\in \omega$, be 
an infinite set of $n$-tuples.     
For every $k$ and every element of $\bar{a}_k$ we fix some 
representative of it in $G^{\omega}$ and think of $\bar{a}_k$ 
as a matrix with $n$ semi-infinite rows, and entries from $G$, 
such that almost all of them are equal to 1.  
Thus we may treat elements $\bar{a}_k$ as semi-infinite sequences 
$\bar{g}_{k,i}$, $i\in\omega$, over the finite alphabet $G^n$ 
such that almost all elements of the sequence are equal to $\bar{1}$. 
Choosing a subset of $\{ \bar{a}_k :k\in\omega\}$ if necessary, 
we may assume that all sequences $\bar{a}_k$ are represented by 
the same set of tuples $\bar{g}_{k,i}$. 
We can also arrange that for every pair $\bar{a}_k$ and $\bar{a}_l$
and any $\bar{g} \in G^n$ the exponent $m$ divides the number 
$|\{i: \bar{g}_{k,i}=\bar{g}\} |-|\{i: \bar{g}_{l,i}=\bar{g}\} |$. 
Let $\bar{h}_0 ,...,\bar{h}_r$ be an enumeration of tuples 
of $G^n$ occurring in all $\bar{a}_k$. 
We may assume that for any $s<t\le r$ and $k\in \omega$ 
the last appearance of $\bar{h}_t$ in $\bar{a}_k$ is after 
the last appearance of $\bar{h}_s$. 
By $l(\bar{a}_k )$ we denote the maximal $i$ for which 
$\bar{g}_{k,i}\neq \bar{1}$. \parskip0pt 

By Lemma \ref{MH}, there are $i<j$ such that $\bar{a}_i$ 
$*$-embeds into $\bar{a}_j$. 
Let $f:\{ 0,...,l(\bar{a}_i )\}$ 
$\rightarrow\{ 0,...,l(\bar{a}_j )\}$ realize this embedding. 
For $s\le r$ let $i_s$ be the greatest $l$ such that 
$\bar{g}_{j,f(l)}=\bar{h}_s$. 
Then let $I_s \subset\omega$ be the set of all $l$ such that 
$\bar{g}_{i,l}= \bar{h}_s$ and $l$ is not in the image of $f$.
We see that for each $s\le r$, the exponent $m$ divides $|I_s|$. 
\parskip0pt 

To define a required embedding $\beta:\Gamma \rightarrow\Gamma$ 
we describe some rules which determine the $\beta$-images 
of elements of $\Gamma$ of the form $(1,...,1,g,1,...,1,...)$ 
where $g$ is the entry with index $l$ (and thus determine $\beta$).  
When $l<l(\bar{a}_i )$ and $\bar{g}_{j,f(l)}$ also appears in 
$\bar{a}_j$ later with a greater index, we define 
the $\beta$-image of $(1,...,1,g,1,...,1,...)$ by the shift 
of $g$ from the index $l$ to $f(l)$.   
In the case when $l\le l(\bar{a}_i )$ has the property that 
$\bar{g}_{j,f(l)}$ does not appear with a greater index in 
$\bar{a}_j$ (thus $f(l)$ is one of the $i_s$-s) we take 
the element as above to 
$(1,...,1,g,1,...,1,g,1,...,1,g,1,...)$, where the last entry 
of $g$ is of the index $f(l)$ and all other appearances of $g$ 
occupy the indexes of the set $I_s$, where $I_s$ is defined 
by $\bar{g}_{j,f(l)}$ as above.     
If $l>l(\bar{a}_i )$, then the $\beta$-image of the element 
above is defined by the shift of $g$ from the index $l$ to 
$l(\bar{a}_j )-l(\bar{a}_i )+l$. 
This construction guarantees that $\beta$ takes $\bar{a}_i$ 
to $\bar{a}_j$. \parskip0pt 

To see that $\beta$ belongs to the closure of the group 
$\Omega$ in the space $\Gamma^{\Gamma}$ take sufficiently 
large $l>l'>l(\bar{a}_j )$ and consider a permutation $\sigma$ 
of $\{ 0,...,l\}$ which extends $f$ and takes every 
$t\in [l(\bar{a}_i )+1,...,l']$ to $t+l(\bar{a}_j )-l(\bar{a}_i )$.  
This permutation naturally extends to the automorphism 
$\hat{\sigma}\in \Omega$ defined as above.  
When we apply $\hat{\sigma}$ together with the product 
$\alpha_{I_0 ,i_0 ,l+1}\cdot \alpha_{I_1 ,i_1 ,l+1}\cdot...\cdot\alpha_{I_{r},i_{r},l+1}$ 
(see Lemma \ref{2type}) we obtain an automorphism of 
$\Gamma$ which coincides with $\beta$ on elements 
represented by sequences which are trivial for indexes 
greater than $l'$. 
This shows that $\beta$ is approximated by automorphisms from $\Omega$. 
\parskip0pt 

It remains to show that $\beta$ preserves the ordering of $\Gamma$. 
Assume $(h_0 ,...,h_t ,...)$ $<(h'_0 ,...,h'_t ,...)$ and $t_0$ 
is the maximal index $t$ such that $h_{t}<h'_{t}$. 
Let $t_1$ be the maximal index where $\beta (h_0 ,...,h_t ,...)$ 
and $\beta (h'_0 ,...,h'_t ,...)$ have distinct entries.  
By the definition of $\beta$, if $l(\bar{a}_i )\le t_0$, then 
$t_1 =t_0 +l(\bar{a}_j )-l(\bar{a}_i )$. 
Since the $t_0$-entries of $(h_0 ,...,h_t ,...)$ and 
$(h'_0 ,...,h'_t ,...)$ coincide with the $t_1$-entries of 
their $\beta$-images respectively, we see that 
$\beta (h_1 ,...,h_t ,...)< \beta (h'_1 ,...,h'_t ,...)$. 
\parskip0pt 

Consider the case when $t_0\le l(\bar{a}_i )$. 
By the definition of $\beta$ the number $t_1$ 
cannot belong to any $I_s$, $s\le r$. 
Thus $t_1 \in Rng(f)$. 
This implies that $f(t_ 0)=t_1$. 
We see that the $t_0$-entries of $(h_0 ,...,h_t ,...)$ and 
$(h'_0 ,...,h'_t ,...)$ coincide with the $t_1$-entries 
of their $\beta$-images respectively and as above we have 
$\beta (h_0 ,...,h_t ,...)<\beta (h'_0 ,...,h'_t ,...)$.         
${\square}$


\begin{thebibliography}{9}
\bibitem{A} Apps, A.B. (1983). 
On $\omega$-categorical class two groups.  
{\em J.Algebra} 82: 516-538. 
\bibitem{AZ} Ahlbtrandt, G., Ziegler, M. (1986).  
Quasi-finitely axiomatizable totally categorical theories. 
{\em Ann. Pure and Appl. Logic} 30: 63 -82.  
\bibitem{AC} Albert, M., Chowdhury, A. (1999). 
The rationals have an AZ-enumeration. {\em J. London Math. Soc. (2)}  
59: 385 - 395. 
\bibitem{CHr} Cherlin, G., Hrushovski, E. (2003).  
{\em Finite Structures with Few Types}.  
Annals of Mathematics Studies, PUP, Princeton. 
\bibitem{CSW} Cherlin, G., Saracino, D., Wood, C. (1993).  
On homogeneous nilpotent groups and rings.  
{\em Proc. Amer. Math. Soc.} 119: 1289 - 1306. 
\bibitem{evans} Evans D. (1994). Examples of $\aleph_0$-categorical 
structures. In: Kaye, R., Macpherson, D., eds. 
{\em Automorphisms of First-Order Structures}. 
Oxford University Press, pp. 33 - 72. 
\bibitem{Hr} Hrushovski, E. (1989). Totally categorical theories. 
{\em Trans. Amer. Math. Soc.} 313: 131 - 159. 
\end{thebibliography}
\end{document}